\documentclass{birkart}

\usepackage{amsmath}
\usepackage{amssymb}
\usepackage{amsfonts}

 \theoremstyle{definition}

 \theoremstyle{remark}

\numberwithin{equation}{section}
\numberwithin{figure}{section}

\newcommand{\beur}{{\mathfrak B}}
\newcommand{\cau}{{\mathfrak C}}
\newcommand{\proj}{{\mathfrak P}}
\newcommand{\mult}{{\mathfrak M}}
\newcommand{\Top}{{\mathfrak T}}
\newcommand{\four}{{\mathfrak F}}

\newcommand{\hispace}{\mathcal H}

\newcommand{\diff}{{\mathrm d}}
\newcommand{\re}{{\mathrm Re}}

\newcommand{\C}{\mathbb C}
\newcommand{\D}{\mathbb D}

\newcommand{\e}{\mathrm e}
\newcommand{\imag}{{\mathrm i}}
\newcommand{\calA}{{\mathcal A}}

\begin{document}
%
\title[On planar Beurling and Fourier transforms]
{On planar Beurling and Fourier transforms}
\author[Hedenmalm]
{H\aa{}kan Hedenmalm}

\address{Hedenmalm: Department of Mathematics\\
The Royal Institute of Technology\\
S -- 100 44 Stockholm\\
SWEDEN}

\email{haakanh@math.kth.se}

\thanks{Research supported by the G\"oran Gustafsson Foundation.}


\subjclass{Primary 32A25, 32A36; Secondary 46E15, 47A15}

\keywords{Beurling transform}


\begin{abstract}
We study the Beurling and Fourier transforms on subspaces of $L^2(\C)$
defined by an invariance property with respect to the root-of-unity group.
This leads to generalizations of these transformations acting unitarily
on weighted $L^2$-spaces over $\C$. 
\end{abstract}

\maketitle








\section{Introduction}

\noindent\bf Beurling and Fourier transforms. \rm
In this note, we shall study certain extensions of the Beurling and Fourier
transforms on the complex plane $\C$. The {\em Fourier transform} of an 
appropriately area-integrable function $f$ is 
$$\four[f](\xi)=\int_{\C}\e^{-2\imag\, \re[z\bar\xi]}\,f(z)\,\diff A(z),
\qquad \xi\in\C,$$
while the {\em Beurling transform} is the singular integral operator
$$\beur_\C[f](z)=\text{pv}\int_\C \frac{f(w)}{(w-z)^2}
\,\diff A(w),\qquad z\in\C;$$
here ``pv'' stands for ``principal value'', and 
$$\diff A(z)=\frac{\diff x\diff y}{\pi},\qquad z=x+\imag y,$$
is normalized area measure. The two transforms are connected via
$$\four\beur_\C[f](\xi)=\frac{\xi}{\bar\xi}\,\four [f](\xi),\qquad \xi\in\C.$$
By the Plancherel identity, $\four$ is a unitary transformation on $L^2(\C)$,
which is supplied with the standard norm
$$\|f\|^2_{L^2(\C)}=\int_\C|f(z)|^2\,\diff A(z).$$
It is clear from this and the above relationship that $\beur_\C$ is unitary
on $L^2(\C)$ as well.
We recall that an operator $T$ acting on a complex Hilbert space $\hispace$ 
is unitary if $T^*T=TT^*=\text{id}$, where $T^*$ is the adjoint and ``id''  
is the identity operator. Expressed differently, that $T$ is unitary means
that $T$ is a surjective isometry.
\medskip

\noindent\bf Root of unity invariance. \rm
For $N=1,2,3,\ldots$, let $\calA_N$ denote the $N$-th roots of unity, that 
is, the collection of all $\alpha\in\C$ with $\alpha^N=1$. For $n=1,\ldots,N$,
we consider the closed subspace $L^2_{n,N}(\C)$ of $L^2(\C)$ consisting of
functions $f$ having the invariance property
\begin{equation}
f(\alpha z)=\alpha^n f(z),\qquad z\in\C,\,\,\,\alpha\in\calA_N.
\label{eq-inv}
\end{equation}
It is easy to see that $f\in L^2_{n,N}(\C)$ if and only $f\in L^2(\C)$ is
of the form
\begin{equation}
f(z)=z^{n}\,g(z^N),\qquad z\in\C,
\label{eq-inv-2}
\end{equation}
where $g$ some other complex-valued function.

We shall study the Beurling and Fourier transforms on the subspaces
$L^2_{n,N}(\C)$. This will be shown to lead to interesting generalizations
of these transforms to the weighted spaces $L^2_\theta(\C)$, with
norm
$$\|f\|^2_{L^2_\theta(\C)}=\int_\C|f(z)|^2\,|z|^{2\theta}\,\diff A(z).$$
Here, $\theta$ is a real parameter. We apply the results obtained regarding
the Beurling transform to conformal mapping, and obtain Grunsky-type 
identities in the spirit of \cite{BH}. The Grunsky-type identity obtained 
here implies the Prawitz inequality (see, e. g., \cite{HS}, \cite{Mil}) 
as a special case. The Grunsky-type inequality that follows from the 
Grunsky-type identity is more or less equivalent to the general Grunsky 
inequality that forms the backdrop to Louis de Branges' work leading
up to the solution of the Bieberbach conjecture (see \cite{dB}).

\section{The Beurling transform}

\noindent\bf The Cauchy transform. \rm 
The {\em Cauchy transform} $\cau_\C$ is the integral transform
$$\cau_\C[f](z)=\int_\C \frac{f(w)}{w-z}\,\diff A(w),$$
defined for appropriately integrable functions. It is related to 
Beurling transform $\beur_\C$ via
$$\beur_\C[f](z)=\partial_z \cau_\C[f](z),$$
where both sides are understood in the sense of distribution theory. Here,
we use the notation
$$\partial_z=\frac12\bigg(\frac{\partial}{\partial x}-\imag
\frac{\partial}{\partial y}\bigg),\quad \bar\partial_z=
\frac12\bigg(\frac{\partial}{\partial x}+\imag
\frac{\partial}{\partial y}\bigg).$$
\medskip

\noindent\bf The Beurling transform and root-of-unity invariance. \rm
Fix an $N=1,2,3,\ldots$ and an $n=1,\ldots,N$. 
We suppose $f\in L^2_{n,N}(\C)$. Then, by the change of variables formula, 
$$\beur_\C[f](z)=
\text{pv}\int_\C \frac{f(w)}{(w-z)^2}
\,\diff A(w)=\text{pv}\int_\C \frac{\alpha^n}{(\alpha w-z)^2}\,f(w)
\,\diff A(w)=\alpha^{n-2}\,\beur_\C[f](\bar\alpha z),\qquad z\in\C,$$  
for $\alpha\in\calA_N$. 
Taking the average over $\calA_N$, we get the identity
$$\beur_\C[f](z)=\frac1N\,\,\text{pv}\int_\C
\sum_{\alpha\in\calA_N}\frac{\alpha^n}{(\alpha w-z)^2}\,f(w)
\,\diff A(w),\qquad z\in\C.$$
\medskip

\noindent\bf A symmetric sum. \rm Next, we study the sum
$$F(z)=\frac1N\sum_{\alpha\in\calA_N}\frac{\alpha^n}{1-\alpha z}.$$
This sum has the symmetry property
$$F(\beta z)=\bar\beta^n \,F(z),\qquad \beta\in \calA_N,$$
which means that $F$ has the form
$$F(z)=z^{N-n}\,G(z^N).$$
The function $G$ then has a simple pole at $1$, and is analytic everywhere
else in the complex plane. Moreover, $F$ vanishes at infinity, so $G$ vanishes
there, too. This leaves us but one possibility, that $G$ has the form
$$G(z)=\frac{C}{1-z},$$
where $C$ is a constant. It is easily established that $C=1$. This leaves us
with
\begin{equation}
F(z)=\frac1N\sum_{\alpha\in\calA_N}\frac{\alpha^n}{1-\alpha z}=
\frac{z^{N-n}}{1-z^N},\qquad z\in\C.
\label{eq-sum}
\end{equation}
As a consequence, we get
$$H(z):=F(z)+zF'(z)=[zF(z)]'=\frac1N\sum_{\alpha\in\calA_N}\frac{\alpha^n}
{(1-\alpha z)^2}=z^{N-n}\bigg\{\frac{N}{(1-z^N)^2}-\frac{n-1}{1-z^N}
\bigg\}.$$ 
This allows us to compute the sum we need:
$$\frac1N\sum_{\alpha\in\calA_N}\frac{\alpha^n}
{(\alpha w-z)^2}=\frac1{z^2}\,H\bigg(\frac {w}{z}\bigg)=z^{n-2}w^{N-n}
\bigg\{\frac{Nz^N}{(z^N-w^N)^2}-\frac{n-1}{z^N-w^N}\bigg\}.$$
For $f\in L^2_{n,N}(\C)$, we thus get the representation
$$\beur_\C[f](z)=z^{n-2}\,\,\text{pv}\int_\C
\bigg\{\frac{Nz^N}{(z^N-w^N)^2}-\frac{n-1}{z^N-w^N}\bigg\}\,w^{N-n}\,f(w)
\,\diff A(w),\qquad z\in\C.$$
Let $f$ and $g$ be connected via (\ref{eq-inv-2}), and implement this 
relationship into the above formula:
\begin{equation}
\beur_\C[f](z)=z^{n-2}\,\,\text{pv}\int_\C
\bigg\{\frac{Nz^N}{(z^N-w^N)^2}-\frac{n-1}{z^N-w^N}\bigg\}\,w^N\,g(w^N)
\,\diff A(w),\qquad z\in\C.
\label{eq-beursymm}
\end{equation}
A similar expression may be found for the Cauchy transform as well:
\begin{equation}
\cau_\C[f](z)=z^{n-N-1}\int_\C
\frac{w^{N}}{w^N-z^N}\,g(w^N)\,\diff A(w),\qquad z\in\C.
\label{eq-cauchysymm}
\end{equation}
\medskip

\noindent\bf The extended Beurling transform. \rm
Let $\Top_\C$ denote the operator
$${\Top}_\C[h](z)=\frac1{z}\,\cau_\C[h](z),$$
and {\em introduce}, for $0\le\theta\le1$, the modified Beurling transform
\begin{equation}
\beur^{\theta}_\C[h](z)=\beur_\C[h](z)+\Top_\C[h](z),
\qquad z\in\C;
\label{def-beurext}
\end{equation}
here, $h$ is assumed to be a nice enough function so that the above Beurling
and Cauchy transforms make sense. It is easy to check that with
$$h(z)=\frac{z\,g(z)}{|z|^{2-2/N}},$$
where $g$ is connected to $f$ via (\ref{eq-inv-2}), we have
$$\beur_\C[f](z)=z^{N+n-2}\,\beur^{(n-1)/N}_\C[h](z^N),
\qquad z\in\C.$$
The fact that $\beur_\C$ is an isometry becomes the norm identity
\begin{equation}
\int_{\C}|h(z)|^2\,|z|^{2\theta}\,\diff A(z)=
\int_\C\big|\beur^{\theta}_\C[h](z)\big|^2\,|z|^{2\theta}\,\diff A(z),
\label{eq-normid}
\end{equation}
where we suppose that $\theta=(n-1)/N$. However, fractions of this type are 
dense in the interval $[0,1]$, so that (\ref{eq-normid}) extends to all 
$\theta$ with $0\le\theta\le1$.
In other words, for $0\le\theta\le1$, {\em the operator $\beur^{\theta}_\C$ 
is unitary on the space} $L^2_\theta(\C)$, which was defined earlier. 
It is known \cite{PV} that $\beur_\C$ is a bounded operator on 
$L^2_\theta(\C)$ for $-1<\theta<1$ (but not for $\theta=\pm1$). 
This means that for $-1<\theta<1$, both terms in (\ref{def-beurext})
are bounded operators on $L_\theta(\C)$. We suspect that the second term in
(\ref{def-beurext}), the operator ${\Top}_\C$, is compact with small 
spectrum. 
 
\medskip

\noindent\bf Extension to real $\theta$. \rm
We first note that $\mult_z$, multiplication by the independent variable,
is an isometric isomorphism $L^2_{\theta+1}(\C)\to L^2_{\theta}(\C)$ for all
real $\theta$. Therefore, for integers $k$ and $0\le\theta\le1$, the operator
$$\beur^{\theta+k}_\C:=\mult_z^{-k}\beur^{\theta}_\C\mult_z^{k}$$ 
is unitary on $L^2_{\theta+k}(\C)$. It supplies an extension
of $\beur^\theta$ to all real $\theta$. Note that since 
$$\frac{1}{(w-z)^2}+\frac{1}{w-z}=\frac{w}{z(w-z)^2},$$
there no disagreement arising from the points $\theta=0$ and $\theta=1$. 
For $-1\le\theta\le0$, $\beur^\theta_\C$ takes the form
$$\beur^\theta_\C=\beur_\C+\theta\,\Top_\C',$$
where
$$\Top'_\C[f](z)=\cau_\C\bigg[\frac{f}{z}\bigg](z).$$
\medskip

\section{Fourier transforms}

\noindent \bf The Fourier transform and root-of-unity invariance. \rm
The Fourier transform of a function in $L^2(\C)$ is given by
$$\four[f](\xi)=\int_{\C}\e^{-2\imag\, \re[z\bar\xi]}\,f(z)\,\diff A(z),
\qquad \xi\in\C.$$
By the Plancherel identity, we have
$$\|\four[f]\|_{L^2(\C)}=\|f\|_{L^2(\C)},\qquad f\in L^2(\C).$$
Now, suppose $f\in L^2_{n,N}(\C)$, so that $f$ has the invariance 
property (\ref{eq-inv}). Then, by the change of variables formula, 
$$\four[f](\xi)=\int_\C \e^{-2\imag\re[z\bar\xi]}\,f(z)\,\diff A(z)=
{\alpha^n}\int_\C \e^{-2\imag\re[\alpha z\bar\xi]}\,f(z)
\,\diff A(z)=\alpha^n\,\four[f](\bar\alpha \xi),\qquad \xi\in\C,$$  
for $\alpha\in\calA_N$. 
Taking the average over $\calA_N$, we get the identity
$$\four[f](\xi)=\int_\C E_{n,N}(z\bar\xi)\,f(z)
\,\diff A(z),\qquad \xi\in\C,$$
where
$$E_{n,N}(z)=\frac1{N}\sum_{\alpha\in\calA_N}\alpha^n\,
\e^{-2\imag\re[\alpha z]}.$$
This sum has the symmetry property
$$E_{n,N}(\beta z)=\bar\beta^n\,E_{n,N}(z),$$
which means that $E_{n,N}(z)$ has the form
$$E_{n,N}(z)=z^{N-n}\,D_{n,N}(z^N).$$
\medskip

\noindent\bf The extended Fourier transform. \rm
We now {\em introduce} the generalized Fourier transform
$$\four_{n,N}[h](\xi)=\frac{|\xi|^{-2(n-1)/N}}{N}
\int_\C D_{n,N}(z\bar\xi)\,h(z)\,\diff A(z).$$
The special case $N=n=1$ gives the standard Fourier transform.
We connect $f$ and $g$ via (\ref{eq-inv-2}).
It is easy to check that with
$$h(z)=\frac{z\,g(z)}{|z|^{2-2/N}},$$
we get that
$$\int_\C|f(z)|^2\,\diff A(z)=
\frac1N\int_\C|h(z)|^2\,|z|^{2(n-1)/N}\,\diff A(z).$$
Moreover, we have
$$\four[f](\xi)=\xi^{n-1}\bar\xi^{N-1}\,\four_{n,N}[h](\xi^N),$$
so that
$$\int_\C|\four[f](\xi)|^2\diff A(\xi)=\frac1N\int_\C
\big|\four_{n,N}[h](\xi)\big|^2\,|\xi|^{2(n-1)/N}\,\diff A(\xi).$$
The Plancherel identity thus states that {\em $\four_{n,N}$ is a unitary
transformation on $L^2_{(n-1)/N}(\C)$}. 
The inverse transformation is quite similar:
$$\four_{n,N}^{-1}[h](z)=(-1)^{N-n}\,\four_{n,N}[h]
\big((-1)^Nz\big).$$
We need to express the function $D_{n,N}$ in a different manner.
Since
$$\e^{-2\imag\,\re[\alpha z]}=\e^{-\imag\,\alpha z}\,
\e^{-\imag\,\bar\alpha\bar z}=\sum_{j,k=0}^{+\infty}
\frac{(-\imag)^{j+k}}{j!k!}\,\alpha^{j-k}\,z^j\bar z^k$$
and
$$\frac1{N}\sum_{\alpha\in\calA_N}\alpha^m=\sum_{l=-\infty}^{+\infty}
\delta_{m,Nl},$$
where delta stands for the Kronecker delta,
we have
\begin{multline*}
E_{n,N}(z)=\frac1{N}\sum_{\alpha\in\calA_N}\alpha^n\e^{-2\imag\,\re[\alpha z]}
=\frac1{N}\sum_{j,k=0}^{+\infty}\frac{(-\imag)^{j+k}}{j!k!}\,z^j\bar z^k
\sum_{\alpha\in\calA_N}\alpha^{j-k+n}\\
=\sum_{l=-\infty}^{+\infty}
\sum_{j,k=0}^{+\infty}\frac{(-\imag)^{j+k}}{j!k!}\,z^j\bar z^k\,
\delta_{j-k+n,Nl}
=\sum_{l=-\infty}^{+\infty}
\sum_{k=0}^{+\infty}\frac{(-1)^k(-\imag)^{-n+Nl}}{(k-n+Nl)!\,k!}
\,z^{k-n+Nl}\,\bar z^{k}\\
=z^{-n}\sum_{l=-\infty}^{+\infty}
\sum_{k=0}^{+\infty}\frac{(-1)^k(-\imag)^{-n+Nl}}{(k-n+Nl)!\,k!}
\,|z|^{2k}\, z^{Nl},
\end{multline*}
with the understanding that 
$$\frac{1}{m!}=0,\qquad m=-1,-2,-3,\ldots.$$
It now follows that
$$D_{n,N}(z)=\sum_{l=-\infty}^{+\infty}
\sum_{k=0}^{+\infty}\frac{(-1)^k(-\imag)^{-n+Nl}}{(k-n+Nl)!\,k!}
\,|z|^{2k/N}\, z^{l-1}.$$
We should mention that as $E_{n,N}(z)$ is bounded by $1$ in modulus,
we have
$$\big|D_{n,N}(z)\big|\le|z|^{n/N-1},\qquad z\in\C.$$
The formula for $D_{n,N}$ allows us to express the modified Fourier transform 
accordingly (provided $f(z)=O(|z|^{-m})$ as $|z|\to+\infty$ for every positive
integer $m$):
$$\four_{n,N}[f](\xi)=\frac{|\xi|^{-2(n-1)/N}}{N}\sum_{l=-\infty}^{+\infty}
\sum_{k=0}^{+\infty}\frac{(-1)^k(-\imag)^{-n+Nl}}{(k-n+Nl)!\,k!}
\,|\xi|^{2k/N}\, \xi^{l-1}\int_{\C}|z|^{2k/N}\, z^{l-1}\,f(z)\,\diff A(z).$$
\medskip

\noindent\bf An application involving the confluent hypergeometric function. 
\rm Next, we consider a function $f$ of the form
$$f(z)=\bar z^{m}|z|^{2\alpha}\,\e^{-\beta|z|^{2/N}},$$
where $m$ is an integer, $\alpha,\beta$ are real with $\beta>0$, and
$$m+2\alpha+\frac{n-1}{N}>-1.$$
If we choose 
$$\alpha=-\frac{n-1-r}{N},$$
where $r$ is an integer, we obtain by calculation that
\begin{multline*}
\four_{n,N}[f](\xi)=
(-\imag)^{-n+N(m+1)}\,\frac{|\xi|^{-2(n-1)/N}\xi^{m}}{N}
\sum_{k=0}^{+\infty}\frac{(-1)^k}{(k-n+N(m+1))!\,k!}
\,|\xi|^{2k/N}\\
\times\int_{\C}|z|^{2(k-n+r+1)/N+2m}\,\e^{-\beta|z|^{2/N}}
\,\diff A(z)\\
=(-\imag)^{-n+N(m+1)}\,
\frac{\xi^{m}|\xi|^{-2(n-1)/N}}{\beta^{r-n+1+N(m+1)}}\,
\sum_{k=0}^{+\infty}\frac{(k+r-n+N(m+1))!}{(k-n+N(m+1))!\,k!}
\,\bigg(-\frac{|\xi|^{2/N}}{\beta}\bigg)^k\\
=(-\imag)^{-n+N(m+1)}\,
\frac{\xi^{m}|\xi|^{-2(n-1)/N}}{\beta^{r-n+1+N(m+1)}}\,
\frac{(N(m+1)+r-n)!}{(N(m+1)-n)!}\\
\times\,\,{}_{1}F_1\bigg(N(m+1)-n+r+1;N(m+1)-n+1;
-\frac{|\xi|^{2/N}}{\beta}\bigg),
\end{multline*}
where ${}_1F_1$ stands for the standard confluent hypergeometric function.
We have from one of the classical identities that
$${}_{1}F_1\bigg(N(m+1)-n+r+1;N(m+1)-n+1;
-\frac{|\xi|^{2/N}}{\beta}\bigg)=e^{-|\xi|^{2/N}/\beta}\,
{}_{1}F_1\bigg(-r;N(m+1)-n+1;\frac{|\xi|^{2/N}}{\beta}\bigg),$$
where the right hand side is easy to compute for positive $r$, as the sum is
then finite.
As a result of the unitarity of $\four_{n,N}$, we find that (with $\beta=1$
and $M=N(m+1)-n+1$)
\begin{equation*}
\int_0^{+\infty}\Big[{}_{1}F_1\big(-r;M;t\big)\Big]^2
\e^{-2t}\,t^{M-1}\diff t
=\frac{(M+2r-1)![(M-1)!]^2}{2^{M+2r}[(M+r-1)!]^2},
\end{equation*} 
which follows from formula 7.622 of \cite{GrR} as a limit case.
\medskip

\section{Applications of Beurling transforms to conformal mapping}

\noindent\bf Transfer to the unit disk. \rm
We need to introduce some general notation. Let $\mult_F$ denote the 
operator of multiplication by the function $F$. 
We also need the Hilbert space $L^2_\theta(X)$ with the norm
$$\|h\|_{L^2_\theta(X)}^2=\int_X|h(z)|^2\,|z|^{2\theta}\,\diff A(z),$$
where $X$ is some Borel measurable subset of $\C$ with positive area.
In the sequel, we fix $\theta$ to the interval $0\le\theta\le1$.
Fix a bounded simply connected domain $\Omega$ in $\C$, which 
contain the origin, and let $\varphi:\D\to\Omega$ denote the conformal mapping
with $\varphi(0)=0$ and $\varphi'(0)>0$. Let $f\in L^2(\Omega)$, and extend 
it to the whole complex plane so that it vanishes on $\C\setminus\Omega$. 
Let $\beur_\Omega[f]$ denote the restriction to $\Omega$ of $\beur_\C[f]$,
and do likewise to define the operators $\cau_\Omega$, $\Top_\Omega$,
$\Top'_\Omega$, $\beur_\Omega^\theta$, as well as $\beur_\Omega^{-\theta}$.  
We introduce transferred operators on spaces over the unit disk in the 
following fashion. First, we suppose $f\in L^2_\theta(\Omega)$. Then the 
associated function
\begin{equation}
g(z)=\bar\varphi'(z)\,\bigg[\frac{\varphi(z)}{z}
\bigg]^\theta\,f\circ\varphi(z),\qquad z\in\D,
\label{eq-f:g}
\end{equation}
belongs to $L^2_\theta(\D)$, with equality of norms:
$$\|g\|_{L^2_\theta(\D)}=\|f\|_{L^2_\theta(\Omega)}.$$
The transferred Cauchy transform is defined as follows:
\begin{equation}
\cau_\varphi^{\theta}[g](z)=\bigg[\frac{\varphi(z)}{z}
\bigg]^\theta\,\cau_\Omega[f](z)=
\int_\D\bigg[\frac{w\,\varphi(z)}{z\,\varphi(w)}\bigg]^{\theta}
\frac{\varphi'(w)}{\varphi(w)-\varphi(z)}\,g(w)\,\diff A(w).
\label{eq-trcau}
\end{equation}
The transferred modified Beurling transform is defined analogously:
\begin{multline*}
\beur_\varphi^\theta[g](z)=\varphi'(z)\bigg[\frac{\varphi(z)}{z}
\bigg]^\theta\beur_\Omega^\theta[f]\circ\varphi(z)
\\
=\varphi'(z)\bigg[\frac{\varphi(z)}{z}
\bigg]^\theta\bigg\{\beur_\Omega[f]\circ\varphi(z)+\frac{\theta}{\varphi(z)}\,
\cau_\Omega[f]\circ\varphi(z)\bigg\}
=\beur_\varphi^{\theta,0}[g](z)+
\theta\,\frac{\varphi'(z)}{\varphi(z)}\,\cau_\varphi^{\theta}[g](z),
\end{multline*}
where
$$\beur_\varphi^{\theta,0}[g](z)=\text{pv}\int_\D
\bigg[\frac{w\,\varphi(z)}{z\,\varphi(w)}\bigg]^{\theta}
\frac{\varphi'(z)\varphi'(w)}{(\varphi(w)-\varphi(z))^2}\,g(w)\,\diff A(w).$$
It is clear that $\beur_\varphi^\theta$ is a norm contraction on 
$L^2_\theta(\D)$. Let $\proj_\theta$ be the integral operator
$$\proj_\theta[f](z)=\int_\D\bigg[\frac1{(1-z\bar w)^2}+
\frac{\theta}{1-z\bar w}\bigg]\,f(w)\,|w|^{2\theta}\diff A(w);$$
it is the orthogonal projection to the subspace of analytic functions in
$L^2_\theta(\D)$. As both $\beur_\varphi^\theta$ and $\proj_\theta$ are
contractions on $L^2_\theta(\D)$, so is their product 
$\proj_\theta\beur_\varphi^\theta$. It remains to represent the operator
$\proj_\theta\beur_\varphi^\theta$ in a reasonable fashion.
The main observation is that
$$\bigg[\frac{w\,\varphi(z)}{z\,\varphi(w)}\bigg]^{\theta}
\frac{\varphi'(z)\varphi'(w)}{(\varphi(w)-\varphi(z))^2}
=\frac{1}{(w-z)^2}-\theta\bigg[\frac{\varphi'(z)}{\varphi(z)}
-\frac{1}{z}\bigg]\frac{1}{w-z}+O(1)$$
near the diagonal $z=w$, so that
\begin{equation}
\bigg[\frac{w\,\varphi(z)}{z\,\varphi(w)}\bigg]^{\theta}
\frac{\varphi'(z)\varphi'(w)}{(\varphi(w)-\varphi(z))^2}
+\theta\,\frac{\varphi'(z)}{\varphi(z)}\,
\bigg[\frac{w\,\varphi(z)}{z\,\varphi(w)}\bigg]^{\theta}
\frac{\varphi'(w)}{\varphi(w)-\varphi(z)}
=\frac{1}{(w-z)^2}+\frac{\theta}{z(w-z)}+O(1),
\label{eq-diagexp}
\end{equation}
again near the diagonal. We observe that in view of (\ref{eq-diagexp}), 
we get the Grunsky-type identity 
\begin{equation}
\proj_\theta\beur_\varphi^\theta=\beur_\varphi^\theta-\beur_\D
+\proj_\theta\beur_\D+\theta\proj_\theta\Top_\D-\theta\Top_\D.
\label{eq-Grunsky}
\end{equation}
To make the involved operators $\proj_\theta\beur_\D$ and 
$\proj_\theta\Top_\D$ appearing in the right hand side of (\ref{eq-Grunsky}) 
more concrete, it is helpful to know that for $\lambda\in\D$,
$$\proj_\theta[f_\lambda](z)=\bar\lambda|\lambda|^{2\theta}
\int_0^1\bigg[\frac1{(1-t\bar\lambda z)^2}+
\frac{\theta}{1-t\bar\lambda z}\bigg]\,t^\theta\diff t,\quad f_\lambda(z)=
\frac{1}{\lambda-z},$$
while
$$\proj_\theta[g_\lambda](z)=-\theta\,\bar\lambda^2|\lambda|^{2\theta-2}
\int_0^1\bigg[\frac1{(1-t\bar\lambda z)^2}+
\frac{\theta}{1-t\bar\lambda z}\bigg]\,t^\theta\diff t,\quad g_\lambda(z)=
\frac{1}{(\lambda-z)^2}.$$
In view of these relations, we quickly verify that
$$\proj_\theta\beur_\D+\theta\proj_\theta\Top_\D=0.$$
The Grunsky-type identity (\ref{eq-Grunsky}) thus simplifies a bit:
\begin{equation}
\proj_\theta\beur_\varphi^\theta=\beur_\varphi^\theta-\beur_\D-\theta\Top_\D.
\label{eq-Grunsky'}
\end{equation} 
The corresponding Grunsky-type inequality reads
\begin{equation}
\big\|\big(\beur_\varphi^\theta-\beur_\D-\theta\Top_\D\big)[f]
\big\|_{L^2_\theta(\D)}\le 
\|f\|_{L^2_\theta(\D)},\qquad f\in L^2_\theta(\D).
\label{eq-Grunsky''}
\end{equation} 
To get a concrete example of how the Grunsky-type inequality works, we
pick
$$f_\lambda(z)=|z|^{-2\theta}\bigg(\frac1{(1-\bar z\lambda)^2}-
\frac{\theta}{1-\bar z\lambda}\bigg),\qquad z\in\D,$$
and compute
\begin{multline*}
\big(\beur_\varphi^\theta-\beur_\D-\theta\Top_\D\big)[f](z)=
\bigg[\frac{\lambda\,\varphi(z)}{z\,\varphi(\lambda)}\bigg]^{\theta}
\frac{\varphi'(z)\varphi'(\lambda)}{(\varphi(\lambda)-\varphi(z))^2}
-\frac{1}{(\lambda-z)^2}\\
+\theta\,\frac{\varphi'(z)}{\varphi(z)}\,
\bigg[\frac{\lambda\,\varphi(z)}{z\,\varphi(\lambda)}\bigg]^{\theta}
\frac{\varphi'(\lambda)}{\varphi(\lambda)-\varphi(z)}
-\frac{\theta}{z(\lambda-z)}.
\end{multline*}
We see that (\ref{eq-Grunsky''}) in this case assumes the form 
($0\le\theta\le1$)
\begin{multline}
\int_\D\bigg|
\bigg[\frac{\lambda\,\varphi(z)}{z\,\varphi(\lambda)}\bigg]^{\theta}
\frac{\varphi'(z)\varphi'(\lambda)}{(\varphi(\lambda)-\varphi(z))^2}
-\frac{1}{(\lambda-z)^2}\\
+\theta\,\frac{\varphi'(z)}{\varphi(z)}\,
\bigg[\frac{\lambda\,\varphi(z)}{z\,\varphi(\lambda)}\bigg]^{\theta}
\frac{\varphi'(\lambda)}{\varphi(\lambda)-\varphi(z)}
-\frac{\theta}{z(\lambda-z)}\bigg|^2|z|^{2\theta}\diff A(z)
\\
\le\int_{\D}|f_\lambda(z)|^2|z|^{2\theta}\diff A(z)=
\int_{\D}\bigg|\frac1{(1-\bar z\lambda)^2}-
\frac{\theta}{1-\bar z\lambda}\bigg|^2|z|^{-2\theta}\diff A(z)
=\frac{1}{(1-|\lambda|^2)^2}-\frac{\theta}{1-|\lambda|^2}.
\label{eq-ex}
\end{multline}
The special case $\lambda=0$ gives us the inequality of Prawitz (see \cite{HS}
and \cite{Mil}; we assume $\varphi'(0)=1$):
\begin{equation*}
\int_\D\bigg|\varphi'(z)\,
\bigg[\frac{\varphi(z)}{z}\bigg]^{\theta-2}
-1\bigg|^2|z|^{2\theta}\diff A(z)
\le\frac{1}{1-\theta}.
\end{equation*}
\medskip

\noindent\bf A dual version. \rm
We carry out the corresponding calculations on the basis of the fact
that $\beur_\C^{-\theta}$ is unitary on $L^2_{-\theta}(\C)$ for 
$0\le\theta\le1$. In analogy with the above treatment, we connect two 
functions $f,g$ via
\begin{equation}
g(z)=\bar\varphi'(z)\,\bigg[\frac{\varphi(z)}{z}
\bigg]^{-\theta}\,f\circ\varphi(z),\qquad z\in\D.
\label{eq-f:g-d}
\end{equation}
Then $f\in L^2_{-\theta}(\Omega)$ if and only if $g\in L^2_{-\theta}(\D)$, 
with equality of norms:
$$\|g\|_{L^2_\theta(\D)}=\|f\|_{L^2_\theta(\Omega)}.$$
The corresponding transferred Beurling transform assumes the form
\begin{multline*}
\beur_\varphi^{-\theta}[g](z)=\varphi'(z)\bigg[\frac{\varphi(z)}{z}
\bigg]^{-\theta}\beur_\Omega^{-\theta}[f]\circ\varphi(z)
\\
=\varphi'(z)\bigg[\frac{\varphi(z)}{z}
\bigg]^{-\theta}
\bigg\{\beur_\Omega[f]\circ\varphi(z)-\theta\,
\cau_\Omega\bigg[\frac{f}{z}\bigg]\circ\varphi(z)\bigg\}
=\beur_\varphi^{-\theta,0}[g](z)-
\theta\,\varphi'(z)\,\cau_\varphi^{-\theta}\bigg[\frac{g}{\varphi}\bigg](z),
\end{multline*}
where $\beur_\varphi^{-\theta,0}$ and $\cau_\varphi^{-\theta}$ are 
as before (just plug in $-\theta$ in place of $\theta$ in the corresponding
formul\ae{}).
It is clear that $\beur_\varphi^{-\theta}$ is a contraction on 
$L^2_{-\theta}(\D)$.

To cut a long story short, the Grunsky-type identity analogous to 
(\ref{eq-Grunsky'}) reads
\begin{equation}
\proj_{-\theta}\beur_\varphi^{-\theta}
=\beur_\varphi^{-\theta}-\beur_\D+\theta\Top'_\D.
\label{eq-Grunsky-d'}
\end{equation}
Let $\bar\proj_{-\theta}^*$ be the operator
$$\bar\proj_{-\theta}^*[g](z)=|z|^{-2\theta}\int_\D
\bigg(\frac{1}{(1-w\bar z)^2}-\frac{\theta}{1-w\bar z}\bigg)\,g(w)\,
\diff A(w);$$
it is a contraction on $L^2_\theta(\D)$, which can be written
$$\bar\proj_{-\theta}^*=\mult_{|z|^{-2\theta}}\bar\proj_{-\theta}
\mult_{|z|^{2\theta}},$$
where $\bar\proj_{-\theta}$ denotes the orthogonal projection onto the
antiholomorphic functions in $L^2_{-\theta}(\D)$. By forming adjoints,
we find that (\ref{eq-Grunsky-d'}) states that
\begin{equation}
\beur_\varphi^\theta\bar\proj_{-\theta}^*
=\beur_\varphi^\theta-\beur_\D-\theta\Top_\D.
\label{eq-Grunsky-d''}
\end{equation}
We now combine (\ref{eq-Grunsky'}) with (\ref{eq-Grunsky-d''}):
\begin{equation}
\beur_\varphi^\theta-\beur_\D-\theta\Top_\D=\proj_\theta\beur_\varphi^\theta=
\beur_\varphi^\theta\bar\proj_{-\theta}^*=\proj_\theta\beur_\varphi^\theta
\bar\proj_{-\theta}^*.
\label{eq-Grunsky-d'''}
\end{equation}
This means that for full mappings $\varphi$, we have equality in the 
Grunsky-type inequality (\ref{eq-Grunsky''}) if and only if $f(z)$ is of the
form $|z|^{-2\theta}$ times an antianalytic function. In particular,
(\ref{eq-ex}) is an equality for full mappings.


\end{document}